\documentclass[a4,12pt]{amsart}
\usepackage[cmtip,arrow]{xy}
\usepackage{pb-diagram,pb-xy,
amsxtra,amsmath,amsthm,amssymb,amscd,ascmac,amsfonts,
multicol,tabularx,latexsym,mathrsfs}
\usepackage{amsfonts}
\theoremstyle{definition}
\newcommand{\Q}{\Bbb Q}
\newcommand{\Z}{\Bbb Z}

\newcommand{\Gal}{\mathrm{Gal}}
\newtheorem{thm}{Theorem}
\newtheorem{lem}{Lemma}
\newtheorem{prop}{Proposition}

\newtheorem{cor}{Corollary}
\newtheorem{prob}{Problem\!}

\newcommand{\coker}{\mathrm{coker}}

\newcommand{\Cl}{\mathrm{Cl}}
\newdimen\minCDarrowwidth
\date{}
\begin{document}
\title[]
{On the abelian groups which occur as Galois cohomology groups of
global unit groups}
\thanks{
\\
This research is partially supported by JSPS, the Grant-in-Aid for
Scientific Research (C) 21540030, and by Waseda University, the Grant for
Special Research Projects 2012A-030.}
\author{Manabu\ OZAKI}
\maketitle
\section{Introduction} 
For any number field $K$ (number field means a finite extension
fields of $\Q$ in what follows), we denote by $E_K$ the unit group of $K$.
If $K/k$ is a Galois extension of number fields, then we define
the Galois cohomology group $\hat{H}^i(K/k,E_K):=
\hat{H}^i(\Gal(K/k),E_K)$ for $i\in \Z$. Here, $\hat{H}^i(G,M)$
stands for the $i$-th Tate cohomology group
for a finite group $G$ and a $G$-module $M$.

We are interested in Galois cohomology groups
$\hat{H}^i(K/k,E_K)$ for many reasons.
For example, when $K/k$ is an {\it unramified} extension,
such cohomology groups are directly related to
the capitulation of ideals, which is one of the major themes of
algebraic number theory;
If $K/k$ is unramified, then we have
\begin{equation}\label{capker}
\begin{aligned}
H^1(K/k,E_K)&\simeq\ker(j_{K/k}:\Cl_k\longrightarrow \Cl_K^{\Gal(K/k)}), \\
H^2(K/k,E_K)&\simeq\coker(j_{K/k}:\Cl_k\longrightarrow \Cl_K^{\Gal(K/k)}),
\end{aligned}
\end{equation}
where $j_{K/k}$ is the natural map from the class group
$\Cl_k$ of $k$ to that $\Cl_K$ of $K$
induced by the inclusion $k\subseteq K$. 

In the present paper, we shall investigate the following problem:
\begin{prob}
For a finite group $G$ and $i\in \Z$, we define
\begin{equation*}
\mathcal{H}^i(G)\!:=\! \{[\hat{H}^i(K/k,E_K)]\,|\,K/k\mbox{\,:\,unramified
$G$-ext'n.\,of number fields}\},
\end{equation*}
where $[A]$ denotes the isomorphism class of $A$ for any abelian group $A$.
What is $\mathcal{H}^i(G)$ ?
\end{prob}
We first note that a simple observation using Dirichlet's unit theorem
and properties of Tate cohomology groups
shows
\[
\mathcal{H}^i(G)\subseteq\mathcal{A}(\# G),
\]
where
$\mathcal{A}(n)$ stands for the set of all the isomorphism classes of the
finite abelian groups killed by $n$ for $n\in\Z$.

The case $i=-1$ and $G=\Z/p$, $p$ being a prime number,
is the most classical situation, which
was first considered by D. Hilbert \cite{Hil}.
He showed $[0]\not\in\mathcal{H}^{-1}(\Z/p)$ in \cite[Satz 92]{Hil}, 
which is equivalent to 
$\ker j_{K/k}\ne 0$ for 
any unramified $\Z/p$-extension $K/k$ of number fields by \eqref{capker}
(Note that $\hat{H}^{-1}(G,M)\simeq\hat{H}^{1}(G,M)$ if $G$ is cyclic).
This latter claim is stated in \cite[Satz 94]{Hil}.

We will reduce our problem to
a problem of group theory.
At present, this reduction has been completely done
only in the case where $G$ is a finite $p$-group, $p$ being a prime number.
Then, in the case where $G$ is a finite $p$-group,
we will determine $\mathcal{H}^i(G)$ for $i=0,1,2,$ and $4$.
We owe determination of $\mathcal{H}^1(G)$ to 
a series of extensive works by K. W. Gruenberg and A. Weiss
\cite{GW2000}, \cite{GW2006}, \cite{GW2011}.
\section{Description of $\hat{H}^i(K/k,E_K)$}
%
%
%
%
%
%
%
We recall the notion of a {\it splitting module} to
describe cohomology groups of unit groups in terms of structures
of certain Galois groups.

Let 
\begin{equation*}
(\varepsilon)
\ \ \ \ \ 
1\longrightarrow A\longrightarrow
\mathcal{G}\longrightarrow G\longrightarrow 1
\end{equation*}
be a group extension of finite groups with abelian kernel $A$.
We denote by $\gamma_{\varepsilon}\in H^2(G,A)$
the cohomology class associated to $(\varepsilon)$.
Define the $G$-module $M_{(\varepsilon)}$ so that
\begin{equation*}
M_{(\varepsilon)}
=A\oplus\bigoplus_{1\ne\tau\in G}
\Z b_\tau
\end{equation*}
as $\Z$-modules.
Here $\{b_\tau\}$ is a free $\Z$-basis, and that
$G$-action on $M_{(\varepsilon)}$
is given by the natural $G$-module structure of $A$
defined by $(\varepsilon)$ and 
\begin{equation*}
\sigma b_\tau=b_{\sigma\tau}-b_{\sigma}
+f(\sigma,\tau)
\end{equation*}
for $\sigma,\,\tau\in G$,
where $f$ is a $2$-cocycle in the cohomology class
$\gamma_{\varepsilon}$
and we set $b_1:=f(1,1)\in A$.
We call $M_{(\varepsilon)}$ {\it the splitting module associated to}
$(\varepsilon)$, and note that 
the $G$-module isomorphism class of $M_{
(\varepsilon)}$
is independent of the choice of a 2-cocycle $f$
of the cohomology class $\gamma_\varepsilon$.

Then we derive from group extension $(\varepsilon)$
the exact sequence of $G$-modules
\begin{equation}\label{epsilon*}
(\varepsilon^*)\ \ \ \ \ 0\longrightarrow A\longrightarrow M_{(\varepsilon)}
\longrightarrow I_G\longrightarrow 0,
\end{equation}
where $I_G$ denotes the augmentation ideal of $\Z[G]$ 
and the map $M_{(\varepsilon)}\longrightarrow I_G$
is given by $b_\sigma\mapsto \sigma-1\ \ (\sigma\in G)$.

Conversely, 
from exact sequence of $G$-modules 
\[
0\longrightarrow A\longrightarrow M\overset{\varphi}{\longrightarrow}
I_G\longrightarrow 0,
\]
we get the group extension (modulo isomorphisms as group extensions of
$G$ by $A$)
\[
(e^{\dagger})\ \ \ \ \ 
1\longrightarrow A\longrightarrow \mathcal{G}\longrightarrow
G\longrightarrow 1
\]
associated to the cohomology class $[f]\in H^2(G,A)$ of the cocycle
$f$ defined by
\[
f(\sigma,\tau):=
\sigma s(\tau-1)-s(\sigma\tau-1)+s(\sigma-1)\in A\ \ \ (\sigma,\tau\in G),
\]
where $s$ is a fixed section of $\varphi$ as $\Z$-modules.
Then we find that
\begin{equation}\label{d-s}
(\varepsilon^{*\dagger})\simeq(\varepsilon),\  \ \ (e^{\dagger *})\simeq (e).
\end{equation}
as group extensions of $G$ by $A$ and $G$-module extensions of $I_G$ by
$A$, respectively. 

Now we will describe the Galois cohomology groups
$\hat{H}^i(K/k,E_K)$ for unramified Galois extensions $K/k$ 
in terms of group extensions naturally arising
from certain towers of unramified Galois extensions. 

Let $K/k$ be an unramified Galois extension of number fields,
and we denote by $H_K$ the maximal unramified abelian extension of $K$.
Put $G:=\Gal(K/k),\ \mathcal{G}:=\Gal(H_K/k)$, and
$A:=\Gal(H_K/K)$.
Then we have the natural group extension 
\[
(\varepsilon_{K/k})\ \ \ \ \ 1\longrightarrow A\longrightarrow
\mathcal{G}\longrightarrow G\longrightarrow 1.
\]

For any finite group $G$,
we define $\mathcal{M}(G)$ to be the class
of all the $G$-modules $M$ fitting into an exact sequence
$0\rightarrow B\rightarrow M\rightarrow I_G\rightarrow 0$
with a finite $G$-module $B$, 
and put  
\[
\frak{X}^i(G):=\{[\hat{H}^{i}(G,M)]|M\in\mathcal{M}(G)\}.
\]
\begin{prop}\label{prop1}
For any unramified Galois extension $K/k$ of number fields 
with $G=\Gal(K/k)$ and $i\in\Z$, we have
\[
\hat{H}^i(K/k, E_K)\simeq \hat{H}^{i-2}(G,M_{(\varepsilon_{K/k})}).
\]
Hence $\mathcal{H}^i(G)\subseteq\frak{X}^{i-2}(G)$ holds.
\end{prop}

{\it Proof.}\ \ \ In fact, this proposition follows from a special case of
Tate sequence given by Ritter-Weiss \cite{RW1996}.
However we will give a direct proof,
which is based on essentially same method as theirs
specialized to the most simple situation, namely,
the Galois extension considered is unramified. 

Let $J_K$, $U_K$, $C_K$, and $\Cl_K$ be the idele group, the unit idele
group, the idele class group,
and the ideal class group of $K$, respectively.
We will observe the following exact commutative diagram:
\begin{equation*}
\begin{diagram}
\dgARROWLENGTH=0.5cm
\dgHORIZPAD=0.5cm
\node[2]{0}
\arrow{s}
\node{0}
\arrow{s}\\
\node[2]{U_K/E_K}\arrow{e,t}{\sim}\arrow{s}\node{H}\arrow{s}\\
\node{0}\arrow{e}\node{C_K}\arrow{s}\arrow{e}
\node{M_{(\overline{\varepsilon})}}\arrow{s}\arrow{e}\node{I_G}
\arrow{e}\arrow{s,=}\node{0}\\
\node{0}\arrow{e}\node{\Cl_K}\arrow{e}\arrow{s}\node{M_{(\varepsilon)}}\arrow{e}
\arrow{s}\node{I_G}\arrow{e}\node{0}\\
\node[2]{0}\node{0},
\end{diagram}
\end{equation*}
where $(\overline{\varepsilon})$ is the group extension
\[
1\longrightarrow C_K
\longrightarrow \overline{\mathcal{G}}
\longrightarrow G
\longrightarrow 1
\]
associated to the fundamental class $c_{K/k}\in H^2(G,C_K)$,
and 
$(\varepsilon)$ is the group extension
\[
1\longrightarrow \Cl_K
\longrightarrow \mathcal{G}
\longrightarrow G
\longrightarrow 1 
\]
associated to the image of $c_{K/k}$ under the natural map
$H^2(G,C_K)\longrightarrow H^2(G,\Cl_K)$.
It follows from Shafarevich's theorem \cite[Chapter 15, Theorem 6]{AT} that
the group extensions $(\varepsilon)$ and $(\varepsilon_{K/k})$
of $G$ are isomorphic via a morphism inducing
the Artin map $\Cl_K\simeq\Gal(H_K/K)$,
hence  
$M_{(\varepsilon)}\simeq M_{(\varepsilon_{K/k})}$
as $G$-modules.
The $G$-module $U_K$ is cohomologically trivial since
$K/k$ is unramified, and
we know $M_{(\overline{\varepsilon})}$
is also cohomologically trivial
(see, for example, \cite[Theorem (3.1.4)]{NSW}).
Therefore we have
\begin{align*}
\hat{H}^i(G,E_K)&\simeq \hat{H}^{i-1}(G, U_K/E_K)
\simeq \hat{H}^{i-1}(G,H)\\
&\simeq \hat{H}^{i-2}(G,M_{(\varepsilon)})\simeq
\hat{H}^{i-2}(G,M_{(\varepsilon_{K/k})})\in\frak{X}^{i-2}(G),
\end{align*}
since $M_{(\varepsilon_{K/k})}\in\mathcal{M}(G)$.
\hfill$\Box$

We make a remark on the case where $G$ is a $p$-group. 
In this case, we naturally define the $p$-quotient of a group
extension
\[
(\varepsilon)\ \ \ 
1\longrightarrow A\longrightarrow
\mathcal{G}\longrightarrow G\longrightarrow 1
\]
with finite abelian $A$ to be
\[
(\varepsilon)_p\ \ \ 
1\longrightarrow A\otimes_\Z\Z_p\longrightarrow
\mathcal{G}_p\longrightarrow G\longrightarrow 1,
\]
where $\mathcal{G}_p$ is the maximal $p$-quotient of $\mathcal{G}$. 
Under this situation, we see that
$\hat{H}^i(G,M_{(\varepsilon)})\simeq \hat{H}^i(G,M_{(\varepsilon)_p})$.
Hence it is enough for our problem
to take account of only the extensions of the forms
\[
1\longrightarrow A\longrightarrow
\mathcal{G}\longrightarrow G\longrightarrow 1
\]
and
\[
0\longrightarrow A\longrightarrow M\longrightarrow I_G\longrightarrow 0
\]
with $\# A$ being a power of $p$ if $G$ is a $p$-group.
\section{Reduction to group theory}
In the case where $G$ is a $p$-group,
we will reduce our problem in Introduction to a problem of group theory
by using the following fact:
\begin{thm}[\cite{Oz}]\label{thm1}
For any finite $p$-group $\mathcal{G}$,
there exists a number field $k$ such that
\[
\Gal(L_p(k)/k)\simeq \mathcal{G}, 
\]
where $L_p(k)$ stands for the maximal unramified $p$-extension
of $k$.
\end{thm}
We can immediately derive the following from the above theorem:
\begin{cor}\label{cor-Oz}
For any given group extension of finite $p$-groups
\[
(\varepsilon)\ \ \ \ \ 1\longrightarrow A\longrightarrow
\mathcal{G}\longrightarrow G\longrightarrow 1
\]
with abelian $A$, 
there exists an unramified $G$-extension $K/k$ such that
$(\varepsilon_{K/k})_p\simeq(\varepsilon)$, namely, there exists
a group isomorphism $\alpha:\mathcal{G}\simeq\Gal(L_p(K)/k)$ 
such that
\begin{equation*}
\begin{CD}
\dgARROWLENGTH=0.5cm
1 @>>> \Gal(L_p(K)/K) @>>> \Gal(L_p(K)/k) @>>> \Gal(K/k) @>>> 1\\
@. @VV{\beta}V @VV{\alpha}V @VV{\gamma}V @.\\
1 @>>>A @>>>\mathcal{G} @>>> G @>>> 1
\end{CD}
\end{equation*}
is an exact commutative diagram,
where $\beta$ and $\gamma$ are the isomorphisms induced by $\alpha$.
\end{cor}
\begin{thm}\label{thm2}
For any finite $p$-group $G$, we have
$\mathcal{H}^i(G)=\frak{X}^{i-2}(G)$ for $i\in\Z$.
\end{thm}
{\it Proof.}\ \ It follows from Proposition 1
that 
$\mathcal{H}^i(G)\subseteq\frak{X}^{i-2}(G)$.

Conversely, let
$[\hat{H}^{i-2}(G,M)]\in\frak{X}^{i-2}(G)\ (M\in\mathcal{M}(G))$
be any element.
Then there exists a group extension
$(\varepsilon)$
such that $M \simeq
M_{(\varepsilon)}$ by \eqref{d-s}.
We derive from Corollary \ref{cor-Oz}
that
there exists an unramified $G$-extension $K/k$ of number fields
such that $(\varepsilon)_p\simeq (\varepsilon_{K/k})_p$.
Hence, by using Proposition \ref{prop1} and the remark after it, we have
\[\ \ \ \ \ 
[\hat{H}^{i-2}(G,M)]=[\hat{H}^{i-2}(G,M_{(\varepsilon_{K/k})})]
=[\hat{H}^{i}(K/k,E_K)]\in \mathcal{H}^{i}(G).\ \ \ \ \ \Box
\]

Thanks to this theorem,
our problem is completely reduced to a purely group theoretic
problem in the case where $G$ is a finite $p$-group.
It is highly interesting and difficult
to search whether $\mathcal{H}^i(G)=\frak{X}^{i-2}(G)$
holds for general finite groups $G$.
\section{Cases $i=2$ and $i=4$.}
Now we will investigate our problem for finite $p$-groups $G$
by using Theorem 2.
We start with rather easy cases:
\begin{thm}
For any finite $p$-group $G$,
we have
\[
\mathcal{H}^2(G)=\mathcal{H}^4(G)=\mathcal{A}(\# G).
\]
\end{thm}
{\it Proof.}\ \ 
We first show that $\mathcal{H}^2(G)=\mathcal{A}(\# G)$.
For any $[X]\in\mathcal{A}(\# G)$, let
$M:=X\oplus I_G\in\mathcal{M}(G)$, 
where we make $G$ act on $X$ trivially.
Then we have $[X]=[\hat{H}^0(G,M)]\in\frak{X}^0(G)=\mathcal{H}^2(G)$
by Theorem \ref{thm2},
because $\hat{H}^0(G,I_G)=0$ and $\hat{H}^0(G,X)\simeq X$.
This shows $\mathcal{A}(\# G)\subseteq\mathcal{H}^2(G)$.
Converse inclusion clearly holds. 

Next we prove $\mathcal{H}^4(G)=\mathcal{A}(\# G)$.
For any $[X]\in\mathcal{A}(\# G)$,
we view $X$ as a $G$-module with trivial $G$-action as above.
We choose a surjection $(\Z/\# G)[G]^{\oplus r}\twoheadrightarrow X$
as $G$-modules
and get an exact sequence
\[
0\longrightarrow Y\longrightarrow (\Z/\# G)[G]^{\oplus r}
\longrightarrow X\longrightarrow 0
\] 
of $G$-modules.
By the same manner, we get a exact sequence
\[
0\longrightarrow Z\longrightarrow (\Z/\# G)[G]^{\oplus s}
\longrightarrow Y\longrightarrow 0
\] 
of $G$-modules for some $s\ge 0$ and finite $G$-module $Z$.
Then we have $Z\oplus I_G\in\mathcal{M}(G)$ and 
\[
\hat{H}^2(G,Z\oplus I_G)\simeq \hat{H}^2(G,Z)\simeq 
\hat{H}^1(G,Y)\simeq
\hat{H}^{0}(G,X)\simeq X
\]
because $H^2(G,I_G)=0$.
Thus we have shown $[X]\in\frak{X}^{2}(G)=\mathcal{H}^4(G)$ by Theorem 2.
This implies our claim. 
\hfill$\Box$
\section{Case $i=1$.}
The case $i=1$ is the most classical situation and many group
theoretic approaches have been available:  

The following theorem is one of the most striking results
given after the principal ideal theorem was proved
by Frutw\"{a}ngler:
\begin{thm}[H.Suzuki \cite{Suz}]\label{thm4}
For any finite abelian group of order $n$, 
we have
\[
\hspace{3.1cm}
\frak{X}^{-1}(G)\subseteq\{[X]\in\mathcal{A}(n)\ |\ n\,|\,\#X\}.
\hspace{3.1cm}\Box
\]
\end{thm}
\begin{cor}\label{cor2}
If $K/k$ is an unramified abelian extension,
then we have
\[
[K:k]\mid \#\ker(\Cl_k\longrightarrow \Cl_K).
\]
\end{cor}
{\it Proof.}\ \ \ 
$[\ker(\Cl_k\longrightarrow \Cl_K)]=[\hat{H}^1(K/k,E_K)]
\in\mathcal{H}^1(G)\subseteq\frak{X}^{-1}(G)
\subseteq\{[X]\in\mathcal{A}(n)\ |\ n\,|\,\#X\}$
by Proposition \ref{prop1} and Theorem \ref{thm4}. 
\hfill$\Box$

\

We note that this corollary implies both of Hilbert's Satz 94
and the principal ideal theorem.
In fact, Suzuki's theorem can be strengthened as follows:
\begin{thm}[(K.W.Gruenberg--A.Weiss \cite{GW2000}]
For any abelian group $G$,
we have
\[\hspace{2.7cm}
\frak{X}^{-1}(G)=\{[X]\in\mathcal{A}(\#G)\ |\ \ \# G\,|\,\# X\}.
\hspace{2.7cm}\Box
\]
\end{thm}
Thanks to this theorem and Theorem \ref{thm2}, we obtain;
\begin{thm}
For any finite abelian $p$-group $G$, we have
\[
\hspace{3cm}
 \mathcal{H}^1(G)= \{[X]\in\mathcal{A}(\# G)\ |\ \ \# G\,|\,\# X\}.
\hspace{3cm}
\Box\]
\end{thm}
\

For any given finite group $G$,
Gruenberg-Weiss \cite{GW2006} showed that
there exists an effectively computable finite subset
$\frak{X}^{-1}_{\mathrm{min}}(G)$ of $\frak{X}^{-1}(G)$ 
such that
\[
\frak{X}^{-1}(G)=
\{[X]\in\mathcal{A}(\# G)\, |\,
\mbox{$[X]$ has a quotient in $\frak{X}^{-1}_{\mathrm{min}}(G)$}\}.
\]
Therefore we can determine $\mathcal{H}^1(G)$ effectively 
when $G$ is a finite $p$-group by using Theorem \ref{thm2}
and the work of Gruenberg and Weiss.
\section{Case $i=0$.}
In this section, we will show;
\begin{thm}\label{thm7}
For any finite $p$-group $G$,
we have
\[
\mathcal{H}^0(G)=\mathcal{A}(\#G).
\]
\end{thm}
We reduce the above theorem as follows:
\begin{lem}\label{lem1}
 If $[0]\in\mathfrak{X}^{-2}(G)$ for a finite $p$-group $G$, then
we have $\mathcal{H}^0(G)=\mathcal{A}(\#G)$.
\end{lem}

{\it Proof.}\ \ \ 
It is sufficient to show $\mathcal{A}(\#G)\subseteq \mathcal{H}^0(G)$
For any $[X]\in\mathcal{A}(\# G)$, 
we view $X$ as a $G$-module with trivial $G$-action
and choose an embedding
\[
X\hookrightarrow (\Z/\#G)^{\oplus r}\simeq((\Z/\#G)[G]^{\oplus r})^G
\subseteq(\Z/\#G)[G]^{\oplus r}
\]
as $G$-modules for some $r\ge 0$.
Then we get the exact sequence
\begin{equation}\label{emb}
0\longrightarrow X\longrightarrow (\Z/\#G)[G]^{\oplus r}
\longrightarrow Y\longrightarrow 0
\end{equation}
of $G$-modules with finite $Y$.

Our assumption implies that there is a exact sequence of $G$-modules
\[
0\longrightarrow Z\longrightarrow M\longrightarrow I_G\longrightarrow 0
\]
such that $\# Z<\infty$ and $\hat{H}^{-2}(G,M)=0$.
Then we have the exact sequence of $G$-modules
\[
0\longrightarrow Z\oplus Y\longrightarrow M\oplus Y\longrightarrow
I_G\longrightarrow 0,
\]
which means $M\oplus Y\in\mathcal{M}(G)$,
and it follows from \eqref{emb} that
\[
\hat{H}^{-2}(G,M\!\oplus\! Y)\simeq
\hat{H}^{-2}(G,Y)\simeq\hat{H}^{-1}(G,X)\simeq X.
\]
This implies $[X]\in\frak{X}^{-2}(G)=\mathcal{H}^0(G)$
by Theorem \ref{thm2}. Thus we have shown
$\mathcal{A}(\#G)\subseteq\mathcal{H}^0(G)$.
\hfill$\Box$

\

To show $[0]\in\mathfrak{X}^{-2}(G)$, we give the following
purely group theoretic proposition:
\begin{prop}\label{prop2}
(a)\ \ For any given finite $p$-group $G$,
there exists a surjective group homomorphism 
$\pi:\overline{\mathcal{G}}\twoheadrightarrow G$
such that
$\overline{\mathcal{G}}$ is a pro-$p$-FAB group
(namely, every open subgroup of $\overline{\mathcal{G}}$
has the finite abelianization) with
$H_2(\overline{\mathcal{G}},\Z_p)=0$.

\noindent
(b)\ \ Let $\pi:\overline{\mathcal{G}}\twoheadrightarrow G$
be any surjective homomorphism with properties
stated in (a), and put $N=\ker\pi$. 
Then, for group extension
\[
(\varepsilon)\ \ \ \ \ 1\longrightarrow N/(N,N)\longrightarrow \overline{\mathcal{G}}/(N,N)
\overset{\overline{\pi}}{\longrightarrow}G\longrightarrow 1,
\]
$\overline{\pi}$ being the map induced by $\pi$,
we have $\hat{H}^{-2}(G,M_{(\varepsilon)})=0$, especially,
$[0]\in\frak{X}^{-2}(G)$.
\end{prop}
Though Proposition 2 seems purely group theoretic,
our proof of it employs largely number theory.

We recall two facts from number theory:
We define $G_{\Q,S}(p)$ to be the Galois group of the maximal $p$-extension
over $\Q$ unramified outside $S$ for any set $S$ of primes of $\Q$.
\begin{thm}
For any finite $p$-group $G$,
there are a finite set $S$ of primes away from $p$ and
a surjection $G_{\Q,S}(p)\twoheadrightarrow G$.
\end{thm}
{\it Proof.}\ \ \ 
This theorem follows from
Shafarevich's theorem (see \cite[Chapter 9, Section 6]{NSW})
for inverse Galois problem on solvable extensions over number fields,
noting that we can make the prime $p$ unramified
in a constructed $G$-extension over $\Q$.
\hfill$\Box$
\begin{thm}(\cite[Theorem 4.9]{Fr})
Let $S$ be a finite set of primes of $\Q$.
We assume that the archimedean prime is contained in $S$
if $p=2$. Then we have $H_2(G_{\Q,S}(p),\Z_p)=0$.
\hfill$\Box$
\end{thm}
Part (a) follows immediately from the above two theorems
because $G_{\Q,S}(p)$ is a pro-$p$-FAB group if $p\not\in S$.

\

To prove part (b), we further recall the following facts from number theory:
\begin{thm}[Folk \cite{F}]\label{folk}
Let $K/k$ be a finite $p$-extension of number fields of finite degree,
and we denote by $H_{K,p}/K$ the maximal unramified abelian $p$-extension.
Then we have $E_k\cap N_{H_{K,p}/k}\left(H_{K,p}^\times\right)
\subseteq N_{K/k}\left(E_K\right)$.
\hfill$\Box$
\end{thm}
\begin{lem}\label{lem2}
Let $L/M$ and $M/k$ be unramified $p$-extensions
of number fields of finite degree such that $L/k$ is normal.
Then we have the commutative diagram
\begin{equation*}
\begin{diagram}
\dgARROWLENGTH=0.5cm
\dgHORIZPAD=0.5cm
\node{E_k/\left(E_k\cap N_{L/k}\left({L}^{\times}\right)\right)}
\arrow{e,J}\arrow{s,r,A}
{\mbox{\small $\mathrm{nat. proj.}$}}
\node{H_2(L/k,\Z_p)}\arrow{s,r}{\mbox{\small $\mathrm{co\mbox{-}inf}$}}\\
\node{E_k/\left(E_k\cap N_{M/k}\left({M}^{\times}\right)\right)}
\arrow{e,J}\node{H_2(M/k,\Z_p),}
\end{diagram}
\end{equation*}
where the right vertical map is the co-inflation map, namely,
the map induced by the natural projection
$\Gal(L/k)\longrightarrow\Gal(M/k)$
and the identity map on $\Z_p$.
\end{lem}
{\it Proof.}\ \ \ 
Since $L/k$ is unramified, the theory of number knot gives a certain
canonical isomorphism
\[
H_2(L/k,\Z_p)\simeq \left(k^\times\cap N_{L/k}(J_L)\right)/N_{L/k}(L^\times)
\]
(see \cite[Section 11.3]{Ta})
and we naturally embed 
$E_k/\left(E_k\cap N_{L/k}\left({L}^{\times}\right)\right)$
into the right hand term of the above isomorphism.
Thus we get the upper horizontal map. 
We also obtain the lower horizontal map similarly. 
The commutativity follows from Horie-Horie \cite[p.618, diagram (3)]{Hor}.
\hfill $\Box$

{\it Proof of Proposition \ref{prop2} (b).}\ \ \ 
Let $\pi:\overline{\mathcal{G}}\twoheadrightarrow G$ and 
$N=\ker\pi$ be as in the statement of Proposition \ref{prop2}.
The fact $H_2(\overline{\mathcal{G}},\Z_p)=0$
implies
that there exists an open normal subgroup 
$H$ of $\overline{\mathcal{G}}$
such that $H\subseteq(N,N)$
and
$
H_2(\overline{\mathcal{G}}/H,\Z_p)
\overset{\mathrm{co\mbox{-}inf}}{\longrightarrow}
H_2(\overline{\mathcal{G}}/(N,N),\Z_p)
$
is the zero map.
Choose $k$ such that there exists an isomorphism
$\delta:\Gal(L_p(k)/k)\simeq\overline{\mathcal{G}}/H$
by using Theorem \ref{thm1}, and
define the tower of number fields $k\subseteq K\subseteq M\subseteq L:=L_p(k)$
so that $\delta$ induces $\Gal(L/M)\simeq (N,N)/H$ and $\Gal(L/K)\simeq N/H$.
%
%
%
Then we see that $\delta$ induces
$\Gal(K/k)\simeq \overline{\mathcal{G}}/N\simeq
G$ and $M=H_{K,p}$, and
by using Lemma \ref{lem2}, we get the commutative diagram
\begin{equation*}
\begin{diagram}
\dgARROWLENGTH=0.5cm
\dgHORIZPAD=0.5cm
\node{E_k/\left(E_k\cap N_{L/k}\left({L}^{\times}\right)\right)}
\arrow{e,J}\arrow{s,r,A}
{\mbox{\small $\mathrm{nat. proj.}$}}
\node{H_2(\overline{\mathcal{G}}/H,\Z_p)}\arrow{s,r}{\mbox{\small $\mathrm{co\mbox{-}inf}$}}\\
\node{E_k/\left(E_k\cap N_{H_{K,p}/k}\left({H^{\times}_{K,p}}\right)\right)}
\arrow{e,J}\node{H_2(\overline{\mathcal{G}}/(N,N),\Z_p).}
\end{diagram}
\end{equation*}
Because the right vertical map in the above diagram is the zero map
from our assumption,
we find that 
$E_k/\left(E_k\cap N_{H_{K,p}/k}
\left({H^{\times}_{K,p}}\right)\right)=0$,
which implies 
$E_k\subseteq N_{K/k}(E_K)$ 
by Theorem \ref{folk}.
Therefore we have 
$\hat{H}^{-2}(G,M_{(\varepsilon)})\simeq\hat{H}^0(K/k,E_K)=0$
by Proposition \ref{prop1},
where
$(\varepsilon)$ is the group extension
$1\rightarrow N/(N,N)\rightarrow
\overline{\mathcal{G}}/(N,N)\overset{\pi}{\rightarrow}G
\rightarrow 1$.
\hfill$\Box$

Thus Theorem \ref{thm7} follows from Lemma \ref{lem1}
and Proposition \ref{prop2}.
\hfill$\Box$

Manabu Ozaki,\par\noindent
Department of Mathematics,\par\noindent
School of Fundamental Science and Engineering,\par\noindent
Waseda University,\par\noindent
Ohkubo 3-4-1, Shinjuku-ku, Tokyo, 169-8555, Japan\par\noindent
e-mail:\ \verb+ozaki@waseda.jp+
\end{document}